\documentclass[11pt]{article}

\usepackage[all]{xy}

\newtheorem{tm}{Theorem}[section]
\newtheorem{lm}[tm]{Lemma}

\newtheorem{rmk}[tm]{Remark}
\newtheorem{cor}[tm]{Corollary}
\newtheorem{ex}[tm]{Example}
\newtheorem{fact}[tm]{Fact}
\newtheorem{??}[tm]{Question}
\newtheorem{defi}[tm]{Definition}

\newtheorem{ass}[tm]{Assumption}

\font\tenmsb=msbm10
\font\sevenmsb=msbm7
\font\fivemsb=msbm5
    
\newfam\msbfam
\textfont\msbfam=\tenmsb
\scriptfont\msbfam=\sevenmsb
\scriptscriptfont\msbfam=\fivemsb
\def\Bbb#1{{\fam\msbfam #1}}

\font\teneufm=eufm10
\font\seveneufm=eufm7
\font\fiveeufm=eufm5
\newfam\eufmfam
\textfont\eufmfam=\teneufm
\scriptfont\eufmfam=\seveneufm
\scriptscriptfont\eufmfam=\fiveeufm
\def\frak#1{{\fam\eufmfam\relax#1}}

\def\lorw{\longrightarrow}
\newcommand\n{\noindent}

\newcommand\ci{\cite}
\newcommand\s{\sigma}
\newcommand\rat{{\Bbb Q}}
\newcommand\comp{{\Bbb C}}

\newcommand\zed{{\Bbb Z}}

\newcommand\pn[1]{{\Bbb P}^{#1}}

\newcommand\blacksquare{{\hspace*{\fill} $\fbox{}$}}

\newcommand{\im}{ \hbox{\rm Im} }

\newcommand{\ke}{ \hbox{\rm Ker} }

\newcommand{\csix}[2]{ {\cal H }^{#1} ( {#2} ) }

\newcommand{\ptd}[1]{ \,^{\frak p}\!\tau_{ \leq {#1} } }

\newcommand{\td}[1]{ \tau_{ \leq {#1} } }

\newcommand{\pc}[2]{ \,^{\frak p}\!{\cal H}^{#1}({#2})   }
\newcommand{\pcs}{ \,^{\frak p}\!{\cal H}   }

\newcommand{\be}{\begin{equation}}
\newcommand{\ee}{\end{equation}}
\newcommand{\la}[1]{\label{#1}}

\title{The Hodge theory of character varieties}
\author{
Mark Andrea A.  de Cataldo\thanks{Partially supported by  N.S.A., N.S.F.
and  Simons' summer research  funds.}\, 
}
\date{}

\begin{document}\maketitle

\begin{abstract}
This is a report  on  joint work
with T. Hausel 
and  L. Migliorini, where we prove, for each of the groups $GL_\comp (2), PGL_\comp (2)$ and $SL_\comp (2)$,
that the non-Abelian Hodge theorem 
identifies
the weight filtration on the cohomology of the character variety with the perverse Leray
filtration on the cohomology of the domain of  the Hitchin map.  We review the decomposition theorem,
N\^go's support theorem, the geometric description of the perverse filtration
and the sub-additivity of the Leray filtration with respect to the cup product. 
\end{abstract}

\tableofcontents

\section{Introduction}
\la{intro}
This is an expanded version of notes from my talk
at the conference ``Classical Algebraic Geometry Today,"
M.S.R.I., Berkeley,  January 25-29, 2009.
The talk consisted of  a report on the  joint work
\ci{dhm} with T. Hausel at Oxford
and  L. Migliorini at Bologna.

Following the recommendation of the editors, I would like this  note
to be accessible to non specialists and to give a small glimpse into an active area of research.
  The reader is referred to the introduction of
\ci{dhm} for more details on what follows.

Let $C$ be a nonsingular complex projective curve.
 We consider
the  following two moduli spaces associated with $C$:
${\mathcal M}:= {\mathcal M}_{Dolbeault}:=$ the moduli space
of  stable holomorphic rank two Higgs bundles on $C$
of degree one (see $\S$\ref{st}) and the character variety
${\mathcal M}':= {\mathcal M}_{Betti}$, i.e. the moduli space 
of irreducible complex dimension two representations of
$\pi_1 (C-p)$ subject to the condition that a loop around
the chosen point 
$p\in C$ is sent to $- {\mbox Id}$.
There is an analogous picture  associated with any complex reductive
Lie group $G$ and the above corresponds to the case $G= GL_{\comp} (2)$.
Our paper \ci{dhm} deals only with the cases when $G= GL_{\comp} (2), PGL_{\comp}(2),
SL_\comp (2)$. Both  ${\mathcal M}$ and ${\mathcal M}'$ are quasi-projective  irreducible and nonsingular
of some even dimension $2d$.
While ${\mathcal M}$ depends on the complex structure of $C$,
${\mathcal M}'$ does not. There is a proper flat and surjective map,
the Hitchin map, $h: {\mathcal M} \lorw \comp^d$ with   general
fibers  Abelian varieties of dimension $d$; in particular, ${\mathcal M}$
is not affine: it contains complete  subvarieties of positive dimension. On the other hand, ${\mathcal M}'$    is easily seen
to be affine (it is a GIT quotient of an affine variety).

The non-Abelian Hodge theorem states that the two  moduli spaces
${\mathcal M}_{Dolbeault}$ and ${\mathcal M}_{Betti}$ are naturally
diffeomorphic,
i.e. that  there is a natural diffeomorphism $\varphi: {\mathcal M} 
\simeq {\mathcal M}'$. Since ${\mathcal M}'$ is affine (resp. Stein) and ${\mathcal M}$
is not affine (resp. not Stein), the map $\varphi$ is not algebraic (resp.
not  holomorphic). Of course, we can still  deduce that $\varphi^*$ is a natural isomorphism
on the singular cohomology groups.

Let us point out that 
the mixed Hodge structure on the cohomology groups $H^j({\mathcal M}, \rat)$
is in fact pure, i.e. every  class has type $(p,q)$ with 
$p+q =j$, or equivalently, every  class has weight $j$.
This follows easily from the fact that, due to the nonsingularity
of ${\mathcal M}$, the weights of $H^j ({\mathcal M},\rat)$ must be $\geq j$.
It remains to show that the weights are also $\leq j$:
the variety 
 ${\mathcal M}$ admits the  fiber $h^{-1} (0)$  of the Hitchin map over the origin $0 \in 
 \comp^d$
 as a deformation retract; it follows  that the restriction map in cohomology,
 $H^j(\mathcal M, \rat) \lorw H^j (h^{-1} (0),\rat)$ is an isomorphism
 of mixed Hodge structures; since the central fiber is compact,
 the weights  of $H^j(h^{-1}(0),\rat)$ are $\leq j$, and we are done.

On the other hand, the mixed Hodge structure on the cohomology 
groups $H^j({\mathcal M}', \rat)$  is known to be
non pure 
(cf. \ci{hauselvillegas}), i.e. there are classes
of degree $j$ but weight $>j$.

It follows that the isomorphism $\varphi^*$ is not compatible
with the two weight filtrations ${\mathcal W}$ on $H^*({\mathcal M},\rat)$
 and ${\mathcal W}'$ on $H^*({\mathcal M}', \rat)$. This fact raises the following question:
{\em if we transplant the weight filtration ${\mathcal W}'$ onto $H^*({\mathcal M},\rat)$
via $\varphi^*$, can we interpret the resulting filtration
on $H^*({\mathcal M}, \rat)$, still called ${\mathcal W}'$, in terms of the
geometry of ${\mathcal M}$}?

The main result in \ci{dhm} is Theorem \ref{dhmmtm} below and it gives
 a positive  answer
to the question raised above. 
In order to state this  answer, we need to introduce one more ingredient and
to  make some trivial
renumerations. (In this paper, we only deal with increasing filtrations.)
The Hitchin map $h: {\mathcal M} \lorw \comp^d$ gives rise
to the perverse
Leray filtration $^{\frak p}\!{\mathcal L}= {^{\frak p}\!{\mathcal L}_h}$ on $H^*({\mathcal M}, \rat)$; this is a suitable
variant of the ordinary Leray filtration for $h$; 
for a geometric description of the perverse Leray  filtration see Theorem \ref{gdpf}.
We renumerate the filtration $^{\frak p}\!{\mathcal L}$ so that $1 \in H^0({\mathcal M},\rat)$
is in place zero (see (\ref{defp})); the resulting renumerated filtration 
 on $H^*({\mathcal M}, \rat)$ is denoted by $P$.
 
All the actual  weights appearing in ${\mathcal W}'$ on $H^*({\mathcal M}',\rat)$
turn out to be multiples of four. We renumerate ${\mathcal W}'$
by setting $W'_k:= {\mathcal W}'_{2k}$.
 
Our answer to the question above is:
 {\em the non-Abelian Hodge theorem isomorphism $\varphi^*$ identifies the weight filtration $W'$ on $H^* ({\mathcal M}', \rat)$
with the perverse Leray filtration $P$ on $H^*({\mathcal M},\rat)$:
\[
P=W'.\]}

The nature of these two filtrations
being very different, we find this coincidence intriguing, but 
at present we cannot explain it
beyond the fact that we can observe it.

The proof of Theorem \ref{dhmmtm} uses a few ideas from 
the topology of algebraic maps. Notably, N\^go's support theorem (\ci{ngo}), the geometric
description of the perverse filtration (\ci{decmigpf}) and the explicit
knowledge of the  cohomology ring $H^*({\mathcal M}_{Betti},\rat)$ (\ci{th}) and of its  mixed Hodge structure
(\ci{hauselvillegas}). 

 One of the crucial ingredients we need
is Theorem \ref{fattob}, which may be of independent interest:
it observes  that Ngo's support theorem  for the Hitchin fibration,
i.e. (\ref{kkll}) below, 
can be refined rather sharply, in the rank two  cases we consider,
as follows:
the intersection complexes appearing in  (\ref{kkll})
are in fact sheaves (up to a dimensional shift).

What follows is a summary of the contents of this paper.
\S\ref{dt} is devoted to stating the decomposition theorem for proper maps
of algebraic varieties and to defining the associated ``supports."
\S\ref{st} states N\^go's support theorem (\ci{ngo}, \S7) and sketches a proof of it
in a special case and 
under a very strong  splitting assumption that does not occur in practice; the purpose here is only to explain the main idea
behind this beautiful result. \S\ref{pflht} is a discussion of the main result
of \ci{decmigpf}, i.e. a description of the perverse filtration in cohomology
with coefficients in a complex via the restriction maps
in cohomology obtained by taking   hyperplane sections.
\S\ref{charvar}  states  the main result in \ci{dhm}, Theorem \ref{dhmmtm},  and discusses
some of the other key ingredients in the proof, notably the use
of the sub-additivity of the ordinary  Leray filtration with respect to cup products. 
Since I could not find
a reference  in the literature for this well-known fact,
I have included  a proof of it  in
 the more technical \S\ref{multler}.
 
 \bigskip
 {\bf Acknowledgments.}  I would like to thank D. Nadler
  for remarks on the support theorem, 
  T. Hausel and L. Migliorini for their  very  helpful comments and
 the referee for
  the excellent suggestions  and remarks.

\subsection{Notation}
\label{not}
We work with sheaves of either Abelian groups, or of  rational vector spaces
over complex algebraic varieties. The survey \ci{decmigbams}
is devoted to the decomposition theorem and contains a more detailed discussion of what follows.

 A sheaf  $F$ on a variety   $Y$  is constructible
if there is a  finite partition $Y= \coprod T_i$ into nonsingular locally closed  irreducible
subvarieties
 that 
 is adapted to $F$, i.e. such that each  $F_{|T_i}$ is  a local system (i.e. a locally constant sheaf)
 on $T_i$.
 A  constructible complex $K$  on a variety $Y$  is a bounded complex of sheaves whose cohomology sheaves ${\cal H}^i(K)$ are  constructible.
We denote by $D_Y$ the corresponding full subcategory
of the derived category of sheaves on $Y$.
If $K \in D_Y$, then
$H^i(Y, K)$ denotes the $i$-th cohomology group of $Y$ with coefficients in $K$. Similarly,
for $H^i_c(Y, K)$.   The complex $K[n]$ has  $i$-th entry   $K^{n+i}$
and differential $d^i_{K[n]} = (-1)^n d_K^{n+i}$.

The standard truncation functors are denoted  by $\td{i}$, the perverse 
(middle perversity) ones by  $\ptd{i}$. 
  The perverse cohomology sheaves are denoted $\pc{i}{K}$, $i \in \zed$.
   We make some use of these notions in  $\S$\ref{multler}.
   Recall that if $K \in D_Y$, then $\pc{i}{K} \neq 0$
  for finitely many values of $i \in \zed$.  In general,
  the collection of perverse cohomology sheaves $\{\! \pc{i}{K}\}_{i\in \zed}$
 does not determine the isomorphism class of  $K$ in $D_Y$;
  e.g. if $j:U \to X$ is the open immersion
  of the complement of a point $p$ in a nonsingular surface $X$,
  then the sheaves $j_! \rat_U$ and $\rat_X \oplus \rat_p$, viewed as complexes
   in $D_X$,
  yield the same collection $\pc{0}{-} = \rat_p$, $\pc{2}{-}= \rat_X [2]$.
  On the other hand, the celebrated decomposition theorem (Theorem \ref{dtpm}
  below) implies
   that
  if $f: X \to Y$ is a proper map of algebraic varieties, with $X$ nonsingular
  for example,
  then the direct image complex $Rf_*\rat_X \simeq \oplus_i
  \pc{i}{Rf_* \rat_X}[-i]$: this implies that the perverse cohomology sheaves
  reconstitute, up to an isomorphism, the  direct image complex; more is true:
  each perverse cohomology sheaf splits further into a direct sum
  of  simple intersection complexes  (cf. (\ref{dtsh})).

 We have the following  subcategories of $D_Y$:
 $D_Y^{\leq 0}:= \{ K \, |\; s.t. \;{\cal H}^i(K) =0, \, \forall i >0\}$ and 
 $^{\frak p}\!D^{\leq 0}_Y:= \{ K \, | \; \dim{\rm supp} \,{\cal H}^i(K) \leq -i, \,
 \forall i \in \zed \}.$
 More generally, a perversity $p$ gives rise
to truncation functors $^p\!\td{i}$, subcategories 
$^{p}\!D^{\leq i}_Y$ and cohomology complexes
$^{p}\! \csix{i}{K}$.

Filtrations on Abelian groups $H$ are assumed to be finite: if the filtration $F_{\bullet}$ on $H$ is increasing,
then $F_i H = 0$ for $i \ll 0$ and $F_i H =H$ for $i \gg 0$;  if $F^{\bullet}$ is decreasing,
then   it is the other way around. 
One can switch type  by setting $F_{i}= F^{-i}$.
For $i \in \zed$, the $i$-th graded
objects  are defined by setting $Gr_i^{F} H:= F_iH/F_{i-1} H$.
The increasing standard filtration ${\cal S}$ on $H^j(Y, K)$ is defined by setting
${\cal S}_i H^j(K):= \im\,\{H^j(Y, \td{i} K) \to H^j(Y, K) \}$.  
Similarly, for $^{\frak p}\!{\cal S}$ and more generally for $^p\!{\cal S}$.
These filtrations are the abutment of corresponding spectral sequences.

Let  $f: X \to Y$ be  a map  of varieties. 
The symbol $f_*$ ($f_!$, resp.)  denotes the 
derived direct image $Rf_*$ (with proper supports $Rf_!$, resp.).
Let $C \in D_X$.
The direct image sheaves are denoted $R^jf_* C$.
We have 
$H^j(X, C) = H^j(Y, f_*C),$ $ H^j_c(X, C) = H^j_c(Y, f_! K)$.

The   Leray filtration  is defined by setting ${\cal L}_i H^j(X, C) := {\cal S}_iH^j(Y, f_*C)$
and it is the abutment 
of the Leray spectral sequence.
 Similarly, for 
$H^j_c(X,C)$. Given a  perversity $p$,  we have
the $p$-Leray spectral sequence abutting to the $p$-Leray filtration  $^p\!{\cal L}$. 
We reserve the terms perverse Leray spectral sequence and perverse
Leray filtration
to the case of middle perversity $p= {\frak p}$.
 
 If $X$ is smooth and $f$ is proper, then we let $Y_{reg}\subseteq Y$ be the   Zariski open 
 set of regular values of $f$  and
 we denote by  $R^i$ the local system  $({R^if_* \rat_X})_{|Y_{reg}}$  on 
 $Y_{reg}$.

\section{The decomposition theorem }
\la{dt}
The purpose of this section is to  state the decomposition theorem \ref{dtpm}  
and to introduce the related  notion of supports.

Let $f: X \to Y$ be a map of varieties.
The Leray spectral sequence  
\[
 E^{pq}_2 = H^p(Y, R^qf_* \rat_X) \Longrightarrow H^{p+q}(X, \rat)
\]
relates the operation of taking cohomology on $Y$ to the same operation on $X$.
If we have $E_2$-degeneration, i.e. $E_2= E_{\infty}$, 
then we have an isomorphism 
\be
\la{ncisp}
H^{j}(X, \rat) \, \simeq \,
\bigoplus_{p+q=j}H^p(Y, R^qf_*\rat).\ee

\begin{ex}
\label{resexe}{\rm
Let $f: X \to Y$ be a resolution of the singularities of
the projective variety 
$Y$. Let us assume, as it is often the case,  that
  the mixed Hodge structure on $H^j(Y, \rat)$ 
is not pure for some $j$.  Then $f^*: H^j(Y,\rat) \to H^*(X,\rat)$
is not injective and    $E_2$-degeneration fails; this is because
injectivity would imply the purity of the mixed Hodge structure
on $H^j(Y,\rat)$.
}\end{ex}
\begin{ex}\la{hopfexe}{\rm
Let $f: (\comp^2 - \{(0,0)\})/\zed \to {\comp}{\pn{1}}$ be a Hopf surface (see \ci{bpvdv}) together with
its natural holomorphic proper submersion onto the projective line.   Since the first  Betti number of the Hopf  surface is one and the one of a fiber is two, $E_2$-degeneration fails.
}
\end{ex}

These examples  show that we cannot expect   $E_2$-degeneration,
neither for  holomorphic proper submersions of compact complex manifolds, nor for
projective   maps of complex  projective  varieties. 
On the other hand,
the following result of P. Deligne \ci{deligne68}  shows 
that $E_2$-degeneration is the norm for proper submersions of complex algebraic varieties.

\begin{tm}
\la{deligne68}
Let $f: X \to Y$ be a smooth proper map of complex algebraic varieties. Then
the Leray spectral sequence for $f$  is $E_2$-degenerate. More precisely,
there is an isomorphism in $D_Y$ 
\[f_* \rat_X \,  \simeq \, \bigoplus_i R^if_*\rat_X [-i].\]
\end{tm}

The decomposition theorem is a far-reaching generalization of Theorem \ref{deligne68} that involves 
intersection cohomology, a notion that  we review briefly next.
A complex algebraic variety $Y$  of dimension $\dim_{\comp} Y=n$ carries intersection cohomology groups  $I\!H^*(Y,\rat)$
and $I\!H_c^*(Y,\rat)$
such that

\smallskip
\n
1.
Poincar\'e duality holds: there is a geometrically defined  perfect   pairing
\[I\!H^{n +j}(Y) \times I\!H^{n-j}_c(Y) \lorw \rat.\]

\n
2.
There is the intersection complex 
 $IC_Y$; it is a constructible  complex of sheaves of rational vector spaces on $Y$
 such that: \[I\!H^j(Y, \rat)= H^{j-n}(Y, IC_Y), \qquad
I\!H^j_c(Y,\rat)= H_c^{j-n}(Y, IC_Y).\]
3.
If $Y$ is nonsingular, then $I\!H^*(Y, \rat) = H^*(Y, \rat)$ and   $IC_Y =\rat_Y[n]$
(complex with the one entry $\rat_Y$  in cohomological degree $-n$).

\n
4.
If $ Y^o$ is a non-empty open subvariety of the nonsingular locus  of $Y$
and $L$ is a local system  on $Y^o$,
then we have  the twisted intersection complex  
$IC_Y(L)$ on $Y$  and the  intersection cohomology groups 
$I\!H^j(Y,L)= {H}^{j-n}(Y, IC_Y(L))$
of $Y$
with coefficients in $L$.

\begin{tm}
\label{dtpm}
{\rm {\bf (Decomposition theorem)}}
{\rm (See \ci{bbd},  Th\'eore\`me 6.2.5.)}

\n
Let $f: X\to  Y$  be a proper map of algebraic varieties. Then
\be
\la{dtsh} f_* IC_X \, \simeq \, \bigoplus_{b\in B} IC_{Z_b}(L_b) [d_b]\ee
for an  uniquely determined finite  collection  $B$  of  triples $(Z_b, L_b, d_b)$ such that
$Z_b\subseteq Y$ 
is a closed irreducible subvariety, 
$ L_b \neq 0$  is 
a simple  local system on some non-empty  and nonsingular Zariski  open  $Z^o_b \subseteq Z_b$
and
 $d_b \in \zed$.
\end{tm}

If, in Theorem \ref{dtpm}, we replace ``simple" with ``semisimple," then we obtain a uniquely determined  collection $B'$
by grouping together the terms with the same cohomological shift $[d_b]$ and the same
irreducible subvariety $Z_b$. 

\begin{defi}
\la{defofsu}
{\rm
The varieties $Z_b\subseteq Y$  are called the
 {\em supports} of the map $f: X \to Y$.
}
\end{defi}

The supports  $Z_b$  are {\em among} the  closed irreducible subvarieties $Z \subseteq Y$ with the property that
\begin{enumerate}
\item\la{11z}
$\exists \; \emptyset \neq \,Z^o \subseteq Z$ over which  all the direct image sheaves
  $R^if_*\rat$ are local systems, and
\item\la{22z}
$Z$ is  maximal with this property.
\end{enumerate}

The following example shows that a support may appear
more than once with distinct cohomological shifts. Of course,
that happens already in the situation of Theorem \ref{deligne68};
the point of the example is that this ``repeated support" may be smaller
than the image $f(X)$.

\begin{ex}\la{morethanonce}{\rm
 Let
$f: X \to Y=\comp^3$  be the blowing-up of a point $o \in \comp^3$;
there is   an isomorphism
\[ f_* \rat_X[3]  \,\simeq \, \rat_Y[3] \oplus \rat_o[1] \oplus \rat_o [-1].\]}
\end{ex}

The next example shows that a variety $Z$, in this case 
$Z=v$, that satisfies conditions \ref{11z} and \ref{22z} above,
may fail to be a support.

\begin{ex}\la{failtosupport}
{\rm
Let  $f: X \to Y$ be the small resolution of the three-dimensional
affine cone $Y \subseteq \comp^4$ over a nonsingular quadric surface $ \pn{1} \times \pn{1} \simeq 
{\mathcal Q}\subseteq \pn{3}$, given by the contraction
to the vertex $v \in Y$  of the
zero section in the total space $X$ of the vector bundle
${\mathcal O}_{\pn{1}} (-1)^2$ . In this case,  we have   
\[Rf_*\rat_X [3] = IC_Y.\]
}\end{ex}

The determination of the supports of a proper map 
 is an important and  difficult problem.

\section{Ng\^o's support theorem}
\label{st}
B.C. Ng\^o has proved (\ci{ngo})  the ``fundamental lemma" in the Langlands program.  
This is a major advance in geometric representation theory,
automorphic representation theory and the arithmetic Langlands program.
See \ci{na}.
One of the  crucial
ingredients of the proof is the support Theorem \ref{nst}, whose proof
applies   the decomposition theorem to the
Hitchin map associated with a reductive group and  a  nonsingular projective  curve.
The support theorem is a rather general result concerning a  certain class of fibrations
with general fibers Abelian varieties and  the Hitchin map is an  important example
of such a fibration.

In our paper 
\ci{dhm} (to which I refer the reader for more context and references),  we deal with the
Hitchin map in the rank two case,  i.e. with the reductive groups
$GL_\comp{2}$, $PGL_\comp (2)$ and $SL_\comp (2)$.
The simpler 
geometry allows us  to refine the conclusion (\ref{kkll})  of the support theorem for the Hitchin map
 in the form of Fact \ref{fattob}, which in turn  we use  in \ci{dhm}
 to prove Theorem \ref{dhmmtm}.

In this section, we discuss the support theorem in the case of $GL_\comp (2)$.
This situation is  too-simple in the context of the fundamental lemma,
but it allows us to concentrate on
the main idea underlying the proof of the support theorem,
i.e.   pursuing the action of Abelian varieties on the fibers of the Hitchin map.
In the context of Ng\^o's work, it is critical to work over  finite fields.
We ignore this important aspect and, for the sake of exposition, we make the   oversimplifying Assumption
\ref{sass} and stick with the situation over $\comp$.

 Let $C$ be a compact  Riemann surface of  genus $g \geq 2$. Let 
${\cal M}$ be   the moduli space of stable
rank $2$   Higgs bundles  on $C$ with determinant  of degree one.
In this context, a point $m \in {\mathcal M}$ parametrizes
a stable pair $(E, \varphi)$, where $E$ is a rank two bundle on $C$ 
with $\deg{(\det{E})} =1$ and $\varphi: E \to E\otimes \omega_C$
(where $\omega_C:= T^*_C$  denotes the canonical bundle of $C$)
is a map of bundles, i.e. a section of $\mbox{End} (E) \otimes \omega_C$.
Stability is a technical condition on the  degrees of the sub-bundles
of $E$ preserved by $\varphi$.
Only the parity of $\deg{(\det{E})}$
counts here: there are only two isomorphism classes of such
moduli spaces; the case of even degree yields a singular moduli space
and we do not say anything new  in that case.

Let $d:= 4g-3$.
The variety ${\cal M}$ is  nonsingular, quasi projective and  of dimension
$2d$.
There is a proper and flat map, called the Hitchin map,  onto affine space
\be
\la{hma}
h: {\cal M}^{2d} \lorw {\Bbb A}^{d} \simeq H^0 (C, \omega_C \oplus \omega_C^{\otimes 2}),\ee
which 
is  a completely integrable system. 

Set-theoretically, the map $h: m = (E,\varphi) \mapsto  (\mbox{trace} (\varphi) , \det{\varphi})$, where the trace and determinant
of the twisted endomorphism $\varphi$ are viewed as sections
of the corresponding powers of $\omega_C$. 

A priori, it is far from clear
that this map is proper. This fact was first noted and proved by Hitchin.
It is a beautiful  fact  (also due to Hitchin) that each nonsingular
fiber    ${\cal M}_a:= h^{-1} (a)$, 
 $ a \in{\Bbb A}^d$,   is     isomorphic to the Jacobian $J(C'_a)$
of what is called the 
spectral curve $C'_a$. This  curve  lives on the  surface given 
by  the total
space of the line bundle $\omega_C$ and it is given set-theoretically
as the double cover  of $C$ given by (and this explains the term ``spectral")
\[
C \ni \{c\} \;  \longleftrightarrow \; \{\mbox{the set of eigenvalues of $\varphi_c$}\}
\in \omega_{C,c}.\]
 The genus $g(C'_a)=d$ by Riemann-Roch and
by the Hurwitz formula.

The singular fibers  of the Hitchin map $h: {\mathcal M} \to {\Bbb A}^d$
are, and this is an euphemism,  difficult to handle.

Let $V \subseteq {\Bbb A}^d$ be the open locus over which
the fibers of $h$ are reduced. The sheaf $R^{2d}_V:= (R^{2d}f_*\rat)_{|V}$
is the $\rat$-linearization of the sheaf of finite sets
given by the sets of 
irreducible components of the  fibers over $V$.  
Let $h_V: {\cal M}_V := h^{-1}(V) \to V$ be 
the restriction of the Hitchin map over $V$.

We can now state Ng\^o's support theorem  in the very  special case at hand.
Roughly speaking, it states that over $V$,  the highest direct image
$R^{2d}_V$ is responsible for all the supports. 
  
\begin{tm}\la{nst}
{\rm {\bf (Ng\^o's support theorem)}} 
A closed and  irreducible subvariety  $Z \subseteq V$  appears   as a support $Z_b$
in the  decomposition theorem  {\rm (\ref{dtsh})}
for $h_V$,
if and only if
there is a dense open subvariety  
$Z^o \subseteq Z$ such that   the restriction  $(R^{2d}_V)_{|Z^o}$   is locally constant
and $Z$ is maximal with this property.
\end{tm}

If we further restrict to the open set $U \subseteq V$ where the fibers are reduced and irreducible,
   then the support theorem has the following striking 
    consequence: the only support  on $U$ 
  is $U$ itself.  The decomposition theorem  (\ref{dtsh}) for  $h_U$  takes then  the following form
   (notation as in \S\ref{not})
   \be\la{kkll}
   {h_U}_* \rat_{ {\cal  M}_U} [2d] \, \simeq 
   \bigoplus_{i=-d}^d IC_U (R^{i+d}) [-i].\ee
The open $U$ is fairly large: its complement  has 
   codimension $\geq 2g-3$.

  \medskip
  The remaining part  of this section is devoted to discussing
  the main idea in the proof of the support theorem.
  
  \smallskip
  There is a group-variety ${\cal P}_V \to V$ over $V$ acting
   on the variety ${\cal M}_V \to V$ over $V$, i.e. a commutative diagram
   \[
   \xymatrix{
   {\cal P}_V \times {\cal M}_V \ar[dr] \ar[rr]^a && {\cal M}_V \ar[dl]\\
   & V &}\]
   satisfying the axioms of an action.
   
   Let us describe this situation over a point $v \in V$.
   The fiber ${\cal M}_v$
   is non-canonically  isomorphic to  a suitable compactification of the identity component ${\cal P}_v$
   of the Picard group ${\rm Pic}(C'_v)$ of the possibly singular
   spectral curve $C'_v$. The variety ${\cal M}_v$ parametrizes
   certain torsion free sheaves  of rank and degree one on $C'_v$.
   The group variety  ${\cal P}_v$ acts on ${\cal M}_v$ via tensor product.
  There is an  exact sequence  
(Chevalley devissage)  of algebraic groups of the indicated dimensions
\be\la{chdev}1\lorw R_v^{\delta_v} \lorw {\mathcal P}_v^d \lorw A_v^{d-\delta_v} \lorw 1
\ee
where   $A_v$ is the Abelian variety given by the  Picard variety  of the normalization 
of the spectral curve $C'_v$ and 
 $R_v$ is an affine  algebraic group.
The  sequence (\ref{chdev}) does not split over the complex numbers,
but it splits over a finite field. It turns out that this is enough
in order to prove the freeness result  on which
the proof of the support theorem rests.
In order to explain the main idea, let us
make the following (over)simplifying assumption.
\begin{ass}
\label{sass}
There is a splitting of  the Chevalley devissage {\rm (\ref{chdev})}.
\end{ass}

A splitting induces an action of $A_v$  on ${\cal M}_v$ with finite stabilizers.
There is the rational homology   algebra $H_*(A_v)$
with product given by the  Pontryagin product $H_i(A_v) \otimes H_j(A_v) \to H_{i+j}(A_v)$
induced by  the cross product, followed by push-forward via the multiplication map in $A_v$.
We have the following standard  (\ci{ngo}, p.134, Proposition 7.4.5) 
\begin{fact}
   \la{finst}
   {\rm }
   Let $A \times T \to T$ be an action of  an Abelian variety  $A$ on a variety $T$
   such that all stabilizers are finite. Then 
$H^*_c(T)$ is a  
free
graded $H_*(A)$-module
for the action of the 
rational    homology algebra $H_*(A)$ on  $H_c^*(T)$.
\end{fact}

Our assumptions imply that
\[
\forall \, v \in V, \qquad 
\mbox{\em $H^*({\cal M}_v)$ is a free graded $H_* (A_v)$-module}.
\]

Let $Z$ be a support appearing in the decomposition theorem (\ref{dtsh}) for $h_V$.
Define a finite set of integers as follows
\[
\mbox{Occ} (Z) := \{ n \in \zed \, | \; \exists b \; {\rm s.t.}  \;
Z_b = Z, \; d_b =-n \} \; \subseteq \; [-d,d].
\]
The integers in ${\rm Occ} (Z)$
  are  in one-to-one correspondence with the summands  (\ref{dtsh}) with support
$Z$. By grouping them, we obtain the graded object
\[
{\frak I}_Z\, := \; \bigoplus_{n \in {\rm Occ} (Z)} IC_{Z}(L^n) [-n].\]

Verdier duality is the generalization  of Poincar\'e duality in the context
of complexes. If we apply this duality to (\ref{dtsh}), we deduce that
 $\mbox{Occ} (Z)$ is symmetric about $0$.

Every intersection complex $IC_Y(L)$ on an irreducible variety $Y$ restricts
to $L[\dim{Y}]$ on a suitable non-empty open subvariety $Y^o \subseteq Y$. It follows
that
there is a non-empty  open subvariety  $V^o \subseteq V$ such that every $IC_Z(L^n)$ restricts
to  $L^n[\dim{Z}]$ on $Z^o:= Z \cap V^o$. Let us consider the restriction of
${\frak I}_Z$ to    $Z^o$:
\[
{\frak L} \, : = \; \bigoplus_{n \in \mbox{Occ} (Z)} L^n[\dim{Z}] [-n].\]
If we set $n^+: = \max{ {\rm Occ}(Z)  }$, then, by the aforementioned
symmetry about the origin, the length
$l ({\frak L}) = 2n^+$.

The decomposition theorem (\ref{dtsh}) over $V^o$ implies that
\[
 \forall \, n \in {\rm Occ}(Z), \qquad
\mbox{\em
$L^n$ is a direct summand of $(R^{2d+n - \dim{Z}}h_* \rat$)}_{|Z^o}.\]

Since the fibers of $h$ have dimension $d$, the higher direct images 
$R^j h_* \rat $ vanish  for every $j >2d$. It follows that
the support theorem  is equivalent to the following 

\medskip
\n
{\bf Claim:}
\[
\mbox{\em $L^{n^+}$  is a direct summand
of $(R^{2d}h_*\rat)_{|Z^o}$.}\]
The Claim is equivalent to  having $n^+ - \dim{Z}=0$, and, again by  the
vanishing for  the direct images $R^j h_*\rat$, this is equivalent to having
$n^+ -\dim{Z} \geq 0$.

Let $z \in Z^o$ be any point.
By adding and subtracting  $\dim{A_z}=d-\delta_z$, we  can re-formulate
the support theorem as follows:
\[  \left[ {\rm codim}  \, Z - \delta_z  \right] + \left[ n^+ - (d-\delta_z)\right]\; \geq\; 0.
\]
It is  thus enough to show
 that each of the two quantities in square brackets is $\geq 0$.

 The first inequality $\left[{\rm codim} \,Z - \delta_z \right] \geq 0$ follows from 
 the deformation theory of Higgs bundles and Riemann-Roch on the curve  $C$.
 This point is standard over the complex numbers.  At present, 
 in positive characteristic it requires 
 the extra freedom of allowing poles of fixed but arbitrary high order. We do not address
 this point here.

\medskip
  Since $l({\frak L}) =2n^+$, in order to prove
  the second inequality, we need to show that
 \[
 l({\frak L}) =  l({\frak L}_z) =2n^+ \geq 2(d - \delta_z). \]
Since $l(H_*(A_z))= 2\dim{A_z} = 2 (d-\delta_z)$,  the inequality
would follow if we
could  prove that:
\[\mbox{
{\em ${\frak L}_z$   is a free graded $H_*(A_z)$-module.}}\]
By virtue of the decomposition theorem,
the graded vector space ${\frak L}_z$ is a graded vector subspace
of $H^*({\cal M}_z)$. This is not enough. We need to make sure
that it is a  free $H_*(A_z)$-submodule.
Once it is known that ${\frak L}_z$ is a submodule, then
 its  freeness is an immediate consequence of standard results
from algebra, notably that a projective module over
the local graded commutative algebra $H_*(A_z)$ is free.
Showing that  ${\frak L}_z$ is  $H_*(A_z)$-stable is
a delicate point, for a priori the contributions
from other supports could enter the picture and spoil it. 
This problem is solved by means of a delicate  specialization argument
which we do not discuss here.

\section{The perverse filtration and the Lefschetz hyperplane theorem}\la{pflht}
Let us review the classical  construction that relates
the Leray filtration on the cohomology of the total 
space a fiber bundle  to the filtration by scheleta
on the base.

Let
 $f: X \to Y$  be  a topological fiber bundle where $Y$ is a 
cell complex of  real dimension $n$.
 Let $Y_*:= \{Y_0 \subseteq \ldots \subseteq Y_k \subseteq \ldots  \subseteq  Y_n =Y\}$
 be the filtration by $k$-scheleta.
 Let  $X_*:= \pi^{-1} (Y_*)$ be the corresponding filtration of the total
 space $X$.

If ${\cal L}$ is the increasing Leray filtration associated with $\pi$, then
we have (see  \ci{spanier}, Ch. 9.4)
\be\la{5647}
{\cal L}_i H^j(X, \zed) = \ke \, \{ H^j (X, \zed) \lorw H^j(X_{j-i-1},\zed)\}.
\ee
 The key fact that one needs (see \ci{decmigpf},  \ci{decsommese})
 is the $\pi$-cellularity of $Y_*$, i.e. the fact that 
\be
\la{cellco}
H^j(Y_p, Y_{p-1}, R^qf_*\zed_X) = 0, \qquad \forall \; j \neq p, \;\;\forall \; q.\ee
This condition is verified  since, for each fixed $p$,  
we are really dealing with   bouquets of $p$-spheres.

This classical result can be viewed as a geometric description of the Leray filtration in the sense
that the subspaces  of the Leray  filtration are exhibited as kernels of restrictions maps
to  the pre-images of scheleta.
The following result of D. Arapura \ci{arapura}  gives a geometric
description of the Leray filtration  for a
projective map of quasi-projective varieties: the important point is that the ``scheleta" can be taken to 
be
algebraic subvarieties!
For generalizations of Arapura's result,
see \ci{decsommese}.
  In what follows, for ease of exposition,
we concentrate on the case when the target is affine.

\begin{tm}
\label{arapura}
{\rm {\bf (Geometric description of the Leray filtration)}}
Let $f:X\to Y$ be a proper map of  algebraic  varieties
with $Y$ affine of dimension $n$.  Then there is a filtration
$Y_*$ of   $\,Y$ by closed algebraic subvarieties  $Y_i$ of dimension $i$ such that
{\rm (\ref{5647})} holds.
\end{tm}

\begin{rmk}
\label{mucrmk}
{\rm
The flag $Y_*$ is constructed inductively as follows. 
Choose a closed embedding $Y \subseteq {\Bbb A}^N$.
Each $Y_i$ is a complete intersection of  $Y$ with $n-i$
sufficiently high degree hypersurfaces in special position.
Here ``special" refers to the fact that in order to achieve
the cellularity condition (\ref{cellco}), we
need to trace, as $p$ decreases,
the $Y_{p-1}$ through the  positive-codimension strata of a partition  of $Y_p$
adapted to the  restricted  sheaves $(R^q f_* \zed_X)_{|Y_p}$.
}
\end{rmk}

Theorem \ref{arapura} affords a simple proof
of the following result of M. Saito (\ci{saito}).
Recall that the integral singular cohomology of
complex  algebraic varieties carries
a canonical and functorial mixed Hodge structure (mHs).

\begin{cor}
\label{mhs}
{\rm {\bf (The Leray filtration is compatible with mHs)}}
Let $f: X \to Y$ be a proper map of algebraic varieties with $Y$ quasi projective.
Then the subspaces of the Leray filtration ${\cal L}$ on $ H^q(X, \zed)$
are mixed Hodge substructures.
\end{cor}

\medskip
Let $K$ be a  constructible complex of sheaves on 
an algebraic variety $Y$. We have 
the   perverse  filtration 
 $^p\!{\cal S}_i H^j(Y,K): = \im   \left\{ H^j \left(Y, \ptd{i} K \right) \to H^j(Y, K) \right\}$.
 Let $f: X \to Y$ be a map of algebraic varieties and let  $C \in D_X$.
We have the perverse Leray filtration $^{\frak p}\!{\cal L}_i$  on $H^j(X,C)$,
i.e the perverse filtration $^{\frak p}\!{\cal S}$   on $H^j(Y, f_*C) = H^j(X, C)$.
Similarly, for  $H^j_c(X, C)$.

\begin{rmk}
\la{esedir}
{\rm
In the situation of  the decomposition theorem  \ref{dtpm}, if we take $X$ to be nonsingular
(if $X$ is singular, then replace cohomology with intersection cohomology in what follows), then
the  subspace
$^{\frak p}\!{\cal L}_i H^j(X, C) \subseteq H^j(X, C)$
is  given by the   images, via the chosen splitting, 
of  the  direct sum  of the $j$-th cohomology groups of the  terms
with $-d_b \leq i$.
The general theory implies that this image is independent
of the  chosen splitting. However, 
different  splittings yield different
embeddings of   each of the direct summands into $H^j(X,C)$.}
\end{rmk}

Let $f: X \to Y$ be a map of varieties where $Y$ is a quasi projective variety.
Let $C \in D_X$ and $K \in D_Y$ (integral coefficients).  
The main result of \ci{decmigpf} is a geometric
description of the perverse and perverse Leray  filtrations.
We state  a significative special case only.

\begin{tm}
\la{gdpf}
{\rm {\bf (Geometric  perverse Leray)}}  
Let $f: X \to Y$ be a map of algebraic varieties with $Y$
affine of dimension $n$.  Then there is a filtration $Y_*$ by closed subvarieties
$Y_i$ of dimension $i$ such that if we take $X_*:= f^{-1} Y_*$, then
\[
^{\frak p}\!{\cal L}_i  \, H^j(X, \zed) := \,^{\frak p}\!{\cal S}_i \, H^j(Y, f_*\zed_X) = \ke \, 
\left\{ H^j(X, \zed) \lorw H^j (X_{n+ j- i-1}, \zed ) \right\}.
\]
\end{tm}

The main difference with respect to Theorem \ref{arapura} is that 
$Y_*$  is obtained by  choosing  general 
vs. special 
hypersurfaces (see Remark \ref{mucrmk}).  This choice is needed
in order to deduce  the perverse analogue of the 
 cellularity condition
 (\ref{cellco}), i.e.
 \[ H^j \left(Y_p, Y_{p-1},   \pcs^q (f_* C\right) = 0,
  \;\; \forall\, j \neq 0, \;\; \forall \; q.\]
These vanishing conditions are verified by a systematic
use of the Lefschetz hyperplane theorem for perverse sheaves.
Unlike  \ci{arapura} and \ci{decsommese}, the proof
for compactly supported cohomology is  completely analogous
to  the one for cohomology.

A second  difference, is that we do not need
the map  $f: X \to Y$ to be proper.  
The choice of general hypersurfaces avoids the usual pitfalls of the failure of the base change theorem
(see \ci{decsommese}).

The  discrepancy  ``$+n$"   between 
(\ref{5647})  for Theorem \ref{arapura} and Theorem \ref{gdpf}
boils down to  the fact that for the affine variety $Y$ of dimension $n$,
the cohomology groups $H^j(Y, F)$ with coefficients in a sheaf (perverse sheaf, resp.) $F$
are non-zero only in the interval $[0,n])$ ($[-n,0]$, resp.).

\medskip
This geometric description of the perverse filtration in terms of the  kernels of restriction maps
to subvarieties is amenable to applications to the mixed Hodge theory of algebraic varieties.
For example, the analogue of Corollary \ref{mhs} holds, with the same proof.
For more applications, see \ci{decpfII}.

\section{Character varieties and the Hitchin fibration: $ P =W'$}\la{charvar}
In this section, I  report on 
  \ci{dhm},
where we  prove Theorem \ref{dhmmtm}.
The main ingredients are the geometric description of the perverse filtration
 in Theorem \ref{gdpf} and    the refinement Theorem \ref{fattob} of  the support theorem
  (\ref{kkll})
 in the case at hand.

We have the
 Hitchin map
(\ref{hma}) for the group  $G= GL_\comp (2)$.
There are analogous maps  
 $\check{h}: \check{\cal M}^{6g-6} \to {\Bbb A}^{3g-3}$ 
 for $G= SL_\comp (2)$  
 and
$\widehat{h}: \widehat{\cal M}^{6g-6}\to {\Bbb A}^{3g-3}$ for
$PGL_\comp (2)$.

Though
these three geometries are
closely related,  this is not the place to  detail the  toing and froing from one group
to another. The main point for this discussion is that we have an
explicit description of the cohomology algebra 
$H^*({\cal M}, \rat)$ in view of the 
  canonical isomorphism
\be\la{hjsa}
H^*({\cal M}, \rat) \, \simeq H^*(\widehat{\cal M}, \rat) \otimes H^*(\mbox{Jac}(C), \rat)
\ee
and of  (\ref{reli}) below. In view of (\ref{hjsa}),
 the key cohomological considerations
towards Theorem \ref{dhmmtm} below
can be made in the $PGL_\comp (2)$ case, for they will
imply easily the ones for $GL_{\comp}(2)$ and, with some extra considerations
which we do not address here, the ones for $SL_{\comp}(2)$.
For simplicity, ignoring some of the subtle differences between
the three groups,
let us work with $\hat{h}: \widehat{\cal M}^{6g-6} \to {\Bbb A}^{3g-3}$.
Though  $\widehat{\cal M}$ is  the quotient of a manifold by the action of a finite group,  for our purposes we can safely
pretend it is a manifold. We set $d:= 3g-3$.

In the context of the non-Abelian Hodge theorem (\ci{simpson}), the quasi projective variety
$\widehat{\cal M}$ is usually denoted $\widehat{\cal M}_D$, where $D$ stands for Dolbeault.
This is to contrast it with the moduli space  $\widehat{\cal M}_B$ 
(Betti) of irreducible $PGL_2(\comp)$
representations of the fundamental group of $C$; this is an affine variety.

The non-Abelian Hodge theorem states that there is a natural diffeomorphism
$\varphi: \widehat{\cal M}_B \simeq \widehat{\cal M}_D$. 
The two varieties are not isomorphic as complex  spaces and, 
a fortiori, neither as algebraic varieties:  the latter  contains
the fibers of the Hitchin map, i.e.
$d$-dimensional Abelian varieties,  while the former is affine.

The diffeomorphism $\varphi$ induces an isomorphism of cohomology rings
$\varphi^*: H^*(\widehat{\cal M}_D, \rat) \simeq  H^*(\widehat{\cal M}_B, \rat)$.
This isomorphism is not compatible with the mixed Hodge structures.
In fact, the mixed Hodge structure on every $H^j(\widehat{\cal M}_D, \rat)$ is 
known to be pure (see  \S\ref{intro}),  while the one on $H^j(\widehat{\cal M}_B, \rat)$ is known to be
not pure
(\ci{hauselvillegas}). 

In particular, the   weight filtrations do not correspond to each other via
$\varphi^*$. 
Our main result in \ci{dhm} can be   stated as follows.

\begin{tm}
\label{dhmmtm}
{\rm \bf{($P = W'$)}} In the cases $G= GL_\comp(2), \, PGL_\comp(2), \, SL_\comp(2)$, 
the non-Abelian Hodge theorem  induces an isomorphism in cohomology
that identifies  the weight filtration for the mixed Hodge structure
on the Betti side with the perverse Leray filtration on the Dolbeault side; 
more precisely,  {\rm (\ref{mopr})} below  holds.
\end{tm}

At present, we do not know what happens if the reductive group $G$  has higher rank.
Moreover,  we do not have
a conceptual explanation
for the so-far mysterious  exchange of structure of Theorem \ref{dhmmtm}.

Our paper
\ci{g=1}  deals with a related moduli space, i.e.  the Hilbert scheme
of $n$ points on the cotangent bundle of an elliptic curve, where a similar
  exchange takes place. 

\medskip
Let us try and describe some of the ideas that play a role
in the proof of Theorem \ref{dhmmtm}.  We refer to \ci{dhm} for details and attributions.

By the work of several people, the cohomology ring and the mixed Hodge structure
of 
$H^*(\widehat{\cal M}_B, \rat)$  are known. 
There are tautological classes:
\[ \alpha \in H^2, \quad \{\psi_i\}_{i=1}^{2g(C)} \in H^3, \quad \beta \in H^4\]
which generate
the cohomology ring. 
With respect to the mixed Hodge structure, these classes 
are  of weight $4$ and pure type $(2,2)$.
Every monomial made of $l$ letters  among these
tautological classes has weight  $4l$ and Hodge type $(2l,2l)$, 
i.e. weights are strictly additive for the cup product. In general, weights 
are only sub-additive. 
There is a   graded $\rat$-algebra  isomorphism
\be\la{reli}H^*(\widehat{\cal M}_B, \rat) \,\simeq\, \frac{\rat \left[ \alpha, \{ \psi_i\}, \beta \right]}{I},\ee
where $I$ is a certain  bihomogeneous ideal with respect to  weight and cohomological degree.
In particular, we have a canonical splitting  for the  increasing weight filtration
${\cal W}'$ on $H^j (\widehat{\cal M}_B, \rat)$ (the trivial weight
filtration ${\mathcal W}$ on the pure  $H^j (\widehat{\cal M}_D, \rat)$ plays no role here):
\[
H^j (\widehat{\cal M}_B, \rat) = \bigoplus_{w\geq 0} H^j_w,
\qquad
{\cal W}'_{w} H^j= \bigoplus_{w' \leq w} H^j_{w'}.\]
The   weights
occur in the interval
$[0, 4 d]$ and  they  are multiples of four, i.e.
${\cal W}'_{4k -i} = {\cal W}'_{4k}$ for every $ 0 \leq i \leq 3$.

By virtue of the decomposition theorem (\ref{dtsh}) and of the fact that the Hitchin map
$\widehat{h}$
is flat of relative dimension $d$, 
the  increasing perverse Leray  filtration  $^{\frak p}\!{\cal L}$  has type $[-d, d]$.

In order to compare  ${\cal W}'$ with $^{\frak p}\!{\cal L}$,
we half the wights, i.e. we set $W'_i := {\cal W}'_{2i}$, and  
we translate  $^{\frak p}\!{\cal L}$, i.e. we set
\begin{equation} \la{defp}
 P:= \,^{\frak p}\!{\mathcal L}(- d).\end{equation}
We still denote
these half-weights by $w$.
We  have that  both $W'$ and $P$ have non-zero graded groups $Gr_i$
only in the interval
$i \in [0, 2d]$.
The two modified filtrations could still be completely unrelated. After all,
they live on  the cohomology of  different algebraic varieties!  The precise formulation
of  Theorem
\ref{dhmmtm}  is 
\be
\label{mopr}
P = W'.\ee

\medskip
Let us  describe our approach to the proof.

\medskip
We introduce the notion of perversity and, ultimately,
we show that the perversity equals the weight.
 We say that $0 \neq u \in H^j(\widehat{\cal M}_D, \rat)$
 has {\em perversity $p= p(u)$} if $u \in {P}_{p} \setminus {P}_{p-1}$. 
 By definition,
 $u=0$ can be given any perversity. Perversities are in the interval $[0, 2d]$.
 We write monomials in the tautological classes as $\alpha^r \beta^s \psi^t$,
 where $\psi^t$ is a short-hand for a product of $t$ classes of type
  $\psi$. 
 Then (\ref{mopr}) can be re-formulated as follows:
 \be\la{wwoo}
p( \alpha^r \beta^s \psi^t)  = w(\alpha^r \beta^s \psi^t) =2(r+s+t).\ee
As it turns out, the harder part is to establish the inequality
\be\la{wwoo2}
p( \alpha^r \beta^s \psi^t)   \leq 2(r+s+t),\ee
for once this is done,  the reverse inequality is proved by   a 
kind of  simple pigeonhole trick.  We thus focus on (\ref{wwoo2}).

Recall that $U \subseteq {\Bbb A}^d$ (see (\ref{kkll}))
  is the dense open set where the fibers of $\hat{h}$ are irreducible. We have the following sharp
  estimate 
  \be\la{shaest}
  \mbox{codim}\, ({\Bbb A}^{d} \setminus U) \geq 2g-3.\ee
  
  For every $0 \leq b \leq d$,
let $\Lambda^b \subseteq {\Bbb A}^{d}$ denote  
a general linear section of dimension $b$. 
We have defined the translated perverse Leray filtration ${P}$  on the cohomology
groups of $\widehat{\cal M}$ for the map $\hat{h}$ that fibers  $\widehat{\cal M}$ 
over  
${\Bbb A}^{d}$. We can do so, in a compatible way,  over $U$ and over the
$\Lambda^b$ so that  the restriction maps respect
the  resulting ${P}$ filtrations.  All these increasing  filtrations start at zero
and perversities are in the interval $[0,2d]$.

The test for perversity Theorem \ref{gdpf}, now reads

\begin{fact}
 \label{factg}
 Let $\Lambda^b \subseteq {\Bbb A}^{d}$ be a general linear subspace of dimension $b$.
 Denote by $\widehat{\cal M}_{\Lambda^b}: = \hat{h}^{-1} ( \Lambda^b)$.
 Then
 \[
 u \in {P}_{j -b-1} H^j(\widehat{\cal M})
 \; \Longleftrightarrow \; u_{|\widehat{\cal M}_{\Lambda^{b}}} =0.\]
 \end{fact}

  We need the following
  strengthening (in the special case we are considering)
  of   the support theorem (\ref{kkll}) over $U$. It is obtained
  by a study of the local monodromy of the family of spectral curves
  around the points of $U$.  Let $j: {\Bbb A}^d_{reg} \to U$ be the open embedding
  of the set or regular values of $\hat{h}$.

\begin{tm}\la{fattob}
The intersection complexes $IC_U(R^i)$ are shifted sheaves and we have
\[
\hat{h}_{U*} \rat \simeq \bigoplus j_*R^i[-i].\]
In particular, the translated perverse Leray filtrations  $P$ coincides with the Leray filtration
$\cal L$
on $H^*(\widehat{\cal M}_U, \rat)$, and
on $H^*(\widehat{\cal M}_{\Lambda^b}, \rat)$  for every $b <2g-3$.
\end{tm}
The last statement is a consequence of (\ref{shaest}): we can trace $\Lambda^b$ inside $U$.

 \medskip
 We can now discuss the scheme of proof for (\ref{wwoo2}).
We start  by establishing  the perversities
of the multiplicative generators, i.e. by proving that 
\[
 p(\alpha) = p(\beta)=p( \psi_i)=2.\]
  By Fact \ref{factg},
 we need to show that $\alpha$ vanishes over the empty set, $\psi_i$ over a point,
 and $\beta$ over a line. The first requirement is of course automatic.
 The second one is a result of M. Thaddeus \ci{thaddeus}. He also proved that $\beta$ vanishes over a point,
 but we need more.
  \begin{fact}
\la{betaline}
The class  $\beta$ is zero over a general line ${\frak l}:= \Lambda^1 \subseteq {\Bbb A}^{3g-3}$.
\end{fact}

\n
 {\em Idea of proof.} 
 By (\ref{shaest}), we can choose
 a general line ${\frak l} = \Lambda^1 \subseteq U$.
 Let  $f: \widehat{\cal M}_{\frak l} := \hat{h}^{-1}({\frak l}) \to {\frak l} $. 
In particular, by abuse of notation, we write
$\widehat{\cal M}_{reg}: = \hat{h}^{-1} ( {\Bbb A}^d_{reg})$, where
${\Bbb A}^d_{reg}\subseteq {\Bbb A}^d$ is the Zariski open and dense
set of regular values of the Hitchin map.

 \n
 By Theorem \ref{fattob} (as it turns out, since we are working over a curve,
 here (\ref{kkll})  is enough to reach the same conclusion)
 we have
\[
 Rf_* \rat \simeq \bigoplus  \left(j_* R^i \right)_{|{\frak l}}[-i]. \]
 In particular, there are no skyscraper summands on ${\frak l}$.
A simple spectral sequence argument over the affine curve ${\frak l}$,  implies 
that the restriction map  $ H^4(\widehat{\cal M}_{\frak l}) \to H^4(\widehat{\cal M}_{{\frak l}_{reg}})$
is injective. (Note that this last conclusion would be clearly false if we had a skyscraper  contribution.)
It is enough to show that $\beta_{|{\widehat{\cal M}_{reg}}}$ is zero.
The class $\beta$ is a multiple of $c_2(\widehat{\cal M})$.
On the other hand, since the Hitchin system
is a completely integrable system over the affine space,
the tangent bundle can be trivialized, in the $C^{\infty}$-sense,  over the open set
of regular point $\widehat{\cal M}_{reg}$  using the Hamiltonian vector fields.  \blacksquare
 
\medskip
Having determined the perversity for the multiplicative generators, we  turn  to (\ref{wwoo2}) which we can 
re-formulate
by saying that  perversities
are sub-additive under cup product.

  In general, I do not know if this is the case: see  the discussion following
  the statement of  Theorem \ref{vvbb} and also  Remark 
\ref{notsha}. On the other hand,
the analogous sub-additivity statement for the Leray filtration $\cal L$
is well-known to hold; see Theorem \ref{vvbb}.

Let us outline our procedure to prove  the sub-additivity of perversity 
in our case.
We want to use the test for perversity Theorem \ref{factg}, for the  monomials in
(\ref{wwoo2}).  First we get rid of $\alpha^r$: in fact,  it is a simple  general  fact
that cupping with a class of degree $i$, raises the perversity by at most $i$.  It follows that
we can concentrate on the case $r=0$.

Here is the outline of the final analysis.
\begin{enumerate}
\item In order to use Theorem \ref{fattob},
we need to  make sure that we can test the monomials over  linear sections
${\Lambda}^b$ which can be traced inside  $U$.
\item
Theorem \ref{fattob}, combined with the sub-additivity of the Leray filtration implies that  we have sub-additivity over $\Lambda^b$.
\item
We deduce
 that the sub-additivity upper bound
 on the perversity  
over $\Lambda^b$, when compared with the cohomological degree
of the monomial, forces the restricted monomial to be zero, i.e.
the monomial passes the test and we are done.
\end{enumerate}

The obstacle in Step 1 is the following:
the dimension $b$ of the testing $\Lambda^b$  increases as 
$s+t$  increases.  On the other hand,  by (\ref{shaest}) we need
$b < 2g-3$. There are plenty of monomials for which $b$ exceeds this bound.
We use the explicit nature of the relations
  $I$ (\ref{reli}) to find an upper bound for $s+t$.
The corresponding upper bound  for $b$ is   $b \leq 2g-3$ (sic!)
and the only class that needs to be tested
on a $\Lambda^{2g-3}$ is $\beta^{g-1}$. This  class turns out to 
require a separate  ad-hoc analysis.
Step 2 requires no further comment. Step 3 is standard as it is 
based on the cohomological dimension
of affine varieties with respect to perverse sheaves.

 \section{Appendix: cup product and Leray filtration}
 \label{multler}
We would like to give a (more or less) self-contained proof of Theorem
\ref{vvbb}, i.e. of the  
fact that the cup product is compatible with the Leray spectral sequence.
We have been unable to locate a suitable reference  in the literature.
As it is clear from our discussion in \S\ref{charvar}, this fact is used  in an essential way
in our proof of Theorem \ref{dhmmtm}

As it turns out, the same proof shows that the cup product is  also compatible
with the $p$-Leray spectral sequence for every non-positive   perversity $p \leq 0$, including ${\frak p}$. 
However, this statement
turns out to be rather weak, unless  we are  in the standard case when $p\equiv 0$. For example, in the case of middle perversity,
 it is  off the mark by
$+d$ with respect to the sub-additivity we need in the proof of Theorem
\ref{dhmmtm}, as it only  implies that  $p(\beta^2) \leq 4+d$,
whereas $p(\beta^2) =4$. Nevertheless,  it seems worthwhile to give a unified proof valid
for every $p \leq 0$.

The statement involves the cup product operation on the cohomology
groups with coefficients in  the direct image  complex. It is thus natural to state 
and prove the compatibility  result for the $p$-standard filtration for arbitrary
complexes on varieties. The compatibility for the $p$-Leray filtration 
is then an immediate consequence. 
We employ  freely the language of derived categories.
We work in the context of  constructible complexes
on algebraic varieties and, just to fix ideas,  with integer coefficients. 
Let us set up the notation necessary to state Theorem \ref{vvbb}.

Let $p: \zed \to \zed$ be any function;
we call it a {\em perversity}.
Given a partition  $X= \coprod S_i$  of a variety $X$ into locally closed 
nonsingular subvarieties $S$
(strata),
we set $p(S): = p(\dim{S})$. By considering all possible partitions of $X$,
this data gives rise to a $t$-structure  on $D_X$ (see \ci{bbd}, p. 56).
The standard $t$-structure corresponds to $p(S)\equiv 0$ and
the middle perversity $t$-structure corresponds to the perversity
${\frak p}$ defined by setting ${\frak p}(S) := - \dim{S}$.

For a given perversity $p$,
the subcategories 
$
^p\!D_X^{\leq i} $ for the corresponding $t$-structure  are  defined  as follows

\[^p\!D_X^{\leq 0} = \left\{ K \in D_X\, | \;
{\cal H}^i (K)_{|S} =0, \; \forall i > p(S)\right\},
\qquad
^p\!D_X^{\leq i}:= \,  ^p\!D^{\leq 0}[-i].
\]
If $p=0$, then $^p\!D_X^{\leq 0}=D_X^{\leq 0}$ is given by the complexes with zero
cohomology sheaves in positive degrees.
If $p = {\frak p}$ is the middle perversity, then  one shows easily
that $^{\frak p}\!D_X^{\leq 0}$ is given 
by those complexes $K$ such that $\dim \mbox{supp}\,{\cal H}^{i}(K) \leq -i$.
By using the truncation functors $^p\!\td{i}$, we can define (see \S\ref{not})
the $p$-standard $^p\!{\cal S}$ and the $p$-Leray
 $^p\!{\cal L}$ filtrations.

Let $K,L \in D_X$.
The tensor product complex  $(K\otimes L, d)$ is defined to be
\be\la{dfgh}
(K\otimes L)^i:= \bigoplus_{a+b=i}K^a \otimes L^b, 
\qquad
d (f_a \otimes g_b) = df \otimes g + (-1)^{a}f \otimes dg.\ee
The left derived tensor product is a bi-functor
$\stackrel{\Bbb L}\otimes: D_X \times D_X \to D_X$  defined
by first taking a flat resolution $L' \to L$  and then by  setting
$K \stackrel{\Bbb L}\otimes L := K \otimes L'$.
If we use field coefficients,
then  the left-derived tensor product
coincides with the ordinary tensor product: $\stackrel{\Bbb L}\otimes  =\otimes$.

Let
\be\la{rbnt}
H^a(X, K) \otimes H^b(X,L) \lorw H^{a+b} (X, K \stackrel{\Bbb L}\otimes L)\ee
be the cup product (\ci{ks}, p. 134).

\medskip
The following establishes that the filtration
in cohomology
associated with  a non-positive perversity  is compatible
with the cup product operation (\ref{rbnt}).

\begin{tm}
\la{vvbb}
Let $p\leq 0$ be a non-positive perversity.
The $p$-standard filtration  and, for a map $f: X \to Y$,  the $p$-Leray   filtration  are compatible with  the
cup product:
\[
{^p\!{\cal S}}_i H^a (X, K) \otimes {^p\!{\cal S}}_j H^b(X, L) \lorw 
 {^p\!{\cal S}}_{i+j} H^{a+b} \left(X, K \stackrel{\Bbb L}\otimes L\right),\]
\[
{^p\!{\cal L}}_i H^a (X, L) \otimes 
{^p\!{\cal L}}_j H^b(X, L) \lorw {^p\!{\cal L}}_{i+j} H^{a+b}
\left(X, K \stackrel{\Bbb L}\otimes L\right).\]
\end{tm}

Theorem \ref{vvbb}  is proved in $\S$\ref{spseqmuly}.
Section \ref{cupcap} shows how the cup product and its variants
for cohomology with compact supports are related to each other;
these variants are  listed in 
(\ref{psvr}). Moreover, if we specialize (\ref{psvr})  to the case of constant coefficients,
and also to the case of the dualizing complex,
then we get  the usual cup products
in cohomology (see  the left-hand-side of (\ref{tanti}))
and the usual cap products involving homology and Borel-Moore homology
(see the right-hand-side of (\ref{tanti})).

\begin{rmk}
\la{moreprod}
{\rm  
The obvious variants of the statement of Theorem \ref{vvbb} hold also for
each of  the variants of the cup product
mentioned above. The same is true
for Theorem \ref{rtmu}, which is merely a souped-up version 
of Theorem \ref{vvbb}.
The reader will have no difficulty repeating,
for each of these variants,  the proof of Theorems \ref{vvbb} and \ref{rtmu}
given in $\S$\ref{spseqmuly}.
}
\end{rmk}

\begin{ex}
\label{notsomucho}
{\rm
We consider  only the two  cases  $p\equiv 0$  and $p = \frak p$;
in the former case we drop the index $p=0$. Let $X$
be a nonsingular variety of dimension $d$.

\begin{enumerate}
\item
Let $K=L=\zed_X$.
Then $1 \in {{\cal S}}_0H^0$ and $1\cup 1=1 \in {\cal S}_0H^0$. 

\item
Let $K=L=\zed_X[d]$.  Then $K\stackrel{\Bbb L}\otimes L= \zed_X[2d]$.
While $1 = 1 \cup 1 \in {^{\frak p}\!{\cal S}}_{-d}H^{-2d}(X,  \zed_X[2d])$,  Theorem
\ref{vvbb} only  predicts  $1 \cup 1 \in {^{\frak p}\!{\cal S}}_0 H^{-2d}(X, \zed_X[2d])$.  

\item
Let $K=L=\zed_p $, where    $p \in X$. We have $1_p =1_p \cup 1_p \in 
{^p\!{\cal S}}_0 H^0 (X, \zed_p)$ and this agrees with
the prediction of  Theorem
\ref{vvbb}.

\item
Let $f:X= Y \times F \to Y$ be the projection, 
and  let $K=L=\rat_X$. 
We have that
\[
{\cal L}_i H^a(X, \rat)= \bigoplus_{i' \leq i}\left(H^{a-i'}(Y, \rat) \otimes_\rat 
H^{i'}(F, \rat)\right).\]
In this case, Theorem \ref{vvbb} is a simple consequence
of the compatibility of the K\"unneth formula with the cup product.

\item
Now let us consider ${^{\frak p}\!{\cal L}}$ for the
same projection map $f: Y \times F \to Y$ as above. We have that $Rf_*\rat \simeq \oplus_{i \geq 0} R^i [-i]$,
where $R^i$ is the constant local system $R^if_*\rat$.
Let us assume that $Y$ is nonsingular of pure dimension $d$. Then
 ${^{\frak p}\!{\cal L}}_i= {\cal L}_{i-d}$, where we use the fact that each  $R^i[d]$ is
 a perverse sheaf due to the nonsingularity of $Y$ (which stems from the one of $X$). We  have that $1 \in {^{\frak p}\!{\cal L}}_d H^0(X, \rat)$.
 On the other hand,    Theorem \ref{vvbb}
predicts only  that $1 = 1 \cup 1 \in {^{\frak p}\!{\cal L}}_{2d}H^0(X, \rat)$. 

\end{enumerate}
}
\end{ex}

These examples, which
as the reader can verify are not an illusion
due to  indexing schemes, show that  Theorem \ref{vvbb} is indeed
sharp.  
However,   its conclusions for ${^{\frak p}\!{\cal S}}$ and ${^{\frak p}\!{\cal L}}$ are 
often  off the mark. 
See  also Remark \ref{notsha}.

\begin{rmk}\la{dunno}
{\rm 
I do not know  an example of a map $f: X \to Y$,
with  $X$ and $Y$ nonsingular, $f$ proper and  flat 
of relative dimension $d$, for which 
the cup product on $H^*(X,\rat)$ does not satisfy 
\begin{equation} \la{try}
{^{\frak p}\!{\cal L}}_i \otimes {^{\frak p}\!{\cal L}}_j \lorw
{^{\frak p}\!{\cal L}}_{i+j-d}\end{equation}
(Theorem \ref{vvbb} predicts that the cup product above lands
in the bigger
${^{\frak p}\!{\cal L}}_{i+j}.$)
In the paper \ci{dhm} we need to establish (\ref{try}) for the Hitchin map.
If (\ref{try}) were true a priori, the proof of the main
result of our paper \ci{dhm} could be somewhat shortened.

\n
Note also  that if the shifted perverse Leray filtration $P:={^{\frak p}\!{\cal L}}(-d)$
(see \ref{defp})
for the Hitchin map $h$ coincided a priori with the ordinary Leray filtration
${\mathcal L}$ of the map $h$,
then (\ref{try}) would follow immediately from the case $p=0$ of Theorem 
\ref{vvbb}. At present, we do not know  if $P= {\mathcal L}$ for the Hitchin map.
In general, i.e. for a map $f$ as above, we have $ {\mathcal L} \subseteq
P$, but the inclusion
can be strict: e.g.
the projection to $\pn{1}$ of the blowing-up of $\pn{2}$ at a point,
where the class of the exceptional divisor is in $P_1$, but it is 
not in ${\mathcal L}_1$.
}
\end{rmk}

\subsection{A simple lemma relating tensor product and truncation}
\la{tpatr}
The key  simple fact behind Theorem \ref{vvbb} in the standard case
when $p \equiv 0$ is that
 if two complexes $K,L \in D^{\leq 0}_X$, i.e. they  have non-zero  cohomology sheaves
in non  negative degrees only, then
the same is true for their derived tensor product. 

 Lemma \ref{bbmmyy} shows that 
the K\"unneth spectral sequence
 implies that the analogous fact  is true
for any  any non-positive perversity $p$.

Let us recall the K\"unneth spectral sequence for the 
derived tensor product of complexes of sheaves. 
Define the Tor-sheaves, a collection of 
 bi-functor with variables sheaves $A$ and $B$, by setting  ${\cal T}or_i (A,B)
:= {\cal H}^{-i}(A \stackrel{\Bbb L}\otimes B)$.  We have
${\cal T}or_0 (A,B) =A \otimes B$ and  ${\cal T}or_i (A,B)=0$ for every $i<0$.
Let $K,L \in D_X$.
We have the K\"unneth spectral sequence (\ci{ega}, III.2., 6.5.4.2,  \ci{verdierthesis}, p.7)
\be\la{kftp}
E_2^{st} = \bigoplus_{a+b=t } {\cal T}or_{-s} (\csix{a}{K}, \csix{b}{L})
\Longrightarrow \csix{s+t}{K \stackrel{\Bbb L}\otimes L}.\ee
This sequence lives in the II-III quadrants, i.e. where $s\leq 0$.
The edge sequence gives  a natural map
\be\la{nlo}
\bigoplus_{a+b=t} \csix{a}{K} \otimes \csix{b}{L}  \lorw 
\csix{t}{K \stackrel{\Bbb L}\otimes L}.
\ee

\begin{lm}
\label{bbmmyy} {\rm {\bf (Tensor product and truncation)}}
Let $p\leq 0$ be any non-positive perversity. Then
\[ \stackrel{\Bbb L}\otimes\,: \;  ^p\!D_X^{\leq i} \times \,  ^p\!D_X^{\leq j} \lorw  \, ^p\!D_X^{\leq  i+j}.
\]
\end{lm}
{\em Proof.} We simplify the notation by dropping the decorations $X$ and
$p$.
Since
\[ 
D^{\leq i} \stackrel{\Bbb L}\otimes D^{\leq j} = D^{\leq 0}[-i] 
\stackrel{\Bbb L}\otimes D^{\leq 0} [-j] = D^{\leq 0} \stackrel{\Bbb L}\otimes D^{\leq 0}  [-i-j],\]
it is enough to show that
\[ 
\stackrel{\Bbb L}\otimes \, : \; D^{\leq 0} \times  D^{\leq 0} \lorw D^{\leq 0}.
\]
We need to verify that the equality 
\be\la{rrttyy}
{\cal H}^{q} (K\stackrel{\Bbb L}\otimes L)_{|S} =0, \qquad  \forall q > p(S),\ee
holds as soon as  the same equality is assumed to hold for $K$ and $L$.

\n
It is enough to prove the  analogous equality
for the Tor-sheaves  on the left-hand-side of   (\ref{kftp}).

\n
Let us note that $p\leq 0$ implies that  if  $L \in D^{\leq 0}$, then
${\cal H}^b(L) =0$ for every $b>0$.




\n
Let  $\s \geq 0$ and  consider
\[
\bigoplus_{a+b =q+\s} {\cal T}or_\s \left({ {\cal H}^a(K)_{|S}, {\cal H}^b(L)_{|S}}
\right) , 
\qquad 
\forall q > p(S).\]
If $a > p(S)$, then ${\cal H}^a(K)_{|S} =0$. 

\n
If $a \leq p(S)$, then $b= q-a+\s >\s \geq 0$, so that
${\cal H}^b(L)=0.$
\blacksquare

\subsection{Spectral sequences and multiplicativity}
\label{spseqmuly}
In this section we prove
Theorem \ref{vvbb}   and we also observe
that it is the reflection at the level of the abutted filtrations of the more general statement Theorem \ref{rtmu}
involving spectral sequences. 

The most efficient formulation is perhaps the one involving the
filtered derived category $D_XF$. We shall quote freely from \ci{illusie}, pp. 285-288. 
To fix ideas, we deal with the cup product in cohomology.  The  formulations
for the  other products  in  \S\ref{cupcap}
are  analogous.

Let $(K, F_1)$ and $(L, F_2)$ be two filtered complexes.
The filtered derived tensor product  $(K \stackrel{\Bbb L}\otimes L, F_{12})$
 is defined as follows.
Let $(L', F_2') \to (L,F_2)$ be a left flat filtered resolution.
Define $K \stackrel{\Bbb L}\otimes L:= K \otimes L'$
and define $F_{12}$ to be the product filtration of $F_1$ and $F_2'$.
We have  natural isomorphisms
\be\la{ksaq}
\bigoplus_{s+s'=\s} \left(Gr_{F_1}^{s} K \stackrel{\Bbb L}\otimes  Gr_{F_2}^{s'} L \right)
\,\stackrel{\simeq}\lorw \,  Gr_{F_{12}^{\s} } 
\left(K \stackrel{\Bbb L}\otimes L\right).
\ee
We have the filtered version of \ci{ks}, p.134,
i.e. a map in $D_{pt}F$
\be\la{ghws}
\left(R\Gamma (X, K), F_1\right) \stackrel{\Bbb L}\otimes 
\left(R\Gamma (X, L), F_2\right) \lorw 
\left(R\Gamma (X,  K \stackrel{\Bbb L}\otimes L), F_{12}\right)
\ee
inducing (cf. \ci{ks}, Ex.I.24.a) the filtered cup product map 
\be\la{ddffgg}
\left(H^a(X, K), F_1\right) \otimes  \left(H^b(X, L), F_2\right)  \lorw
\left(H^{a+b} \left(X,  K \stackrel{\Bbb L}\otimes L\right), F_{12}\right).\ee

In view of Theorem \ref{rtmu}, by  first 
recalling the notion of bilinear pairing of spectral sequences
(\ci{spanier}, p. 491),
we have a bilinear pairing of spectral sequences
\be\la{cmdk}
E_1^{st}(K, F_1) \otimes E_1^{s't'}(L, F_2) \lorw E_1^{s+s', t+t'} \left(
K \stackrel{\Bbb L}\otimes L, F_{12}\right)
\ee 
that  on the $E_1$-term  coincides with the cup  product map
(\ref{rbnt})
induced by (\ref{ksaq}),  and on the $E_{\infty}$-term
is the graded cup product 
associated with the filtered cup product  (\ref{ddffgg}).

Given   $(M,F) \in D_XF$, we have the spectral sequence
\[
E_1^{st} = E_1^{st}(M,F)= H^{s+t}(X, Gr_{F}^s M) \Longrightarrow
H^{s+t}(X, M), \qquad
E_{\infty}^{st} = Gr_F^s H^{s+t} (X, M),
\]
with abutment the filtration induced by $(M,F)$ on  the cohomology groups
$H^*(X, M)$.
Clearly, we can always compose with the map of spectral sequences
induced by any filtered map
$(K \stackrel{\Bbb L}\otimes L, F_{12}) \to (M,F)$.

\medskip
We apply the machinery above to the case when the filtrations
$F_i$ are the $p$-standard decreasing 
 filtrations  $^p\!{\cal S}$ induced by the $t$-structure associated with
 a non-positive perversity  $p\leq 0$. 
 The construction of $^p\!{\cal S}$ is performed via the use 
 of injective resolutions (\ci{decmigpf}, \S3.1).
The product filtration $F_{12}$ on the derived tensor product is {\em not}
 the  $p$-standard filtration, not even up to isomorphism in
 the filtered derived category $D_XF$; see Remark \ref{notsha} below.
 
The upshot of this discussion is   that Lemma \ref{bbmmyy} implies the following
\begin{lm}
\la{rtvb}
There is a canonical lift 
 \[u\, :\; \left(K \stackrel{\Bbb L}\otimes L, F_{12}\right) \lorw
 \left(K \stackrel{\Bbb L}\otimes L, \,^p\!{\cal S} \right)\]
 of the identity on $K \stackrel{\Bbb L}\otimes L$ to  $D_XF$.
 \end{lm} 
 {\em Proof.} Let $N$ denote the derived tensor product of $K$ with $L$.
  It is enough to show that $F^{\s}_{12}N \in \, ^p\!D^{\leq  -\s}_X$, for every $\s \in \zed$.
 We prove this by decreasing induction on $\s$. 
 The statement is clearly true for $\s \gg 0$. We have the  short exact sequence
 \[
 0\lorw F^{\s+1}_{12} N \to F^{\s}_{12} N  \to Gr^{\s}_{F_{12}} N \lorw 0.\] 
 Lemma \ref{bbmmyy} implies that  $Gr^{\s}_{F_{12}} N \in  \, ^p\!D^{\leq  -\s}_X$ and the inductive hypothesis gives    $F^{\s+1}_{12} N  \in \,  ^p\!D^{\leq  -\s -1}_X \subseteq 
 \,  ^p\!D^{\leq  -\s}_X$.  We have the following simple fact:
 if $A \to B \to C \to A[1]$ is a distinguished triangle and $A, C \in \, ^p\!D^{\leq  i}_X$,
then $B \in \, ^p\!D^{\leq  i}_X$. We conclude the proof
by applying  this fact to the distinguished triangle
associated with  the short exact sequence above.
 \blacksquare

 \medskip
 We are now ready for the
 
 \medskip
 \n
 {\bf Proof of Theorem \ref{vvbb}:}
 
 \n
 {Apply the construction (\ref{ddffgg}) to the case $F_i= {^p\!{\mathcal S}}$.
 Compose the resulting filtered cup product map with the canonical lift
 $u$ of Lemma \ref{rtvb} and obtain the filtered
 cup product map of Theorem \ref{vvbb}}.\blacksquare

 \medskip
 As mentioned earlier, Theorem \ref{vvbb} is the abutted reflection
 of the following  statement concerning spectral sequences:
 
 \begin{tm}
 \label{rtmu}
 There is a natural  bilinear pairing of spectral sequences
\[
E_1^{st}(K, \, ^p\!{\cal S}) \otimes E_1^{s't'}(L, \,^p\!{\cal S}) \lorw E_1^{s+s', t+t'} 
\left(
K \stackrel{\Bbb L}\otimes L, \, ^p\!{\cal S}\right)
\]
such that:  

\begin{enumerate}
\item
 on the $E_1$-term   it coincides with the cup  product map
induced by {\rm (\ref{ksaq})}, which in this case reads 
\be\la{sldk}
^p\!\csix{-s}{K} [s] \stackrel{\Bbb L}\otimes \, 
 ^p\!\csix{-s'}{K} [s'] \lorw  \, ^p\!{\cal H}^{-s-s'}\left( K \stackrel{\Bbb L}\otimes L\right) [s+s'],
\ee
 \item \la{234}  on the $E_{\infty}$-term it 
is the graded cup product 
associated with  the filtered cup product  {\rm (\ref{ddffgg})}.
\end{enumerate}
\end{tm}
{\em Proof:} 
Compose (\ref{cmdk}) with the map of spectral sequences induced by the canonical map $u$
of Lemma \ref{rtvb}.
 \blacksquare

 \begin{rmk}
 \la{notsha}
 {\rm
 Unless we are in the case $p\equiv 0$, the product  filtration
 $F_{12}$  of the $p$-standard filtrations is often strictly smaller than the $p$-standard
 filtration. As a result, 
 the graded pairing is  often trivial. One can see this  on the $E_1$-page
 in terms of the map (\ref{sldk}).  Here is an example.
 Let $X$ be nonsingular of pure dimension
 $d$, take middle perversity and  perverse complexes $K=L= \rat_X[d]$.
 The pairing in question is
 $\rat_X [d] \otimes \rat_X [d] \lorw \pc{0}{\rat_X} [2d] =0.$
 }
 \end{rmk}

\subsection{Cup and cap}
\label{cupcap} The methods employed in the previous sections are of course
susceptible of being applied to the other usual constructions, such as the cup product in cohomology
with compact supports and cap products in homology and in Borel-Moore homology.  

By taking various flavors of (\ref{ghws}) with compact supports,
we obtain the commutative diagram of cup product maps
\be\la{psvr}
\xymatrix{
H_c^i (X, K) \otimes H_c^j(X, L) \ar[r] \ar[d] \ar[ddr] & H_c^{i+j}(X,K \stackrel{\Bbb L}\otimes L) \ar[d]^= \\
H^i(X,K) \otimes H_c^j(X,L)\ar[r] \ar[d] \ar[dr] & H_c^{i+j}(X, K \stackrel{\Bbb L}\otimes L)\ar[d] \\
H^i (X,K) \otimes H^j(X,L)  \ar[r]  & H^{i+j}(X,K \stackrel{\Bbb L}\otimes L).
}
\ee
Theorem \ref{rtmu} applies to each row, each vertical arrow is a filtered map for the product filtrations
and, as a result, the conclusion of Theorem \ref{rtmu} apply to the diagonal products
as well.

We have the following important special cases.
Take
$K$ and $L$ to be either $\zed_X$ and/or $\omega_X$ (the Verdier dualizing complex of $X$).
We have $\zed_X \stackrel{\Bbb L}\otimes \omega_X =\omega_X$ as well as 
 the following equalities (decorations omitted)
\[
H^i(X, \zed)= H^i(X, \zed_X)=H^i\zed_X= H^i, \quad 
 H^i_c = H^i_c \zed_X, \quad H_i =H^{-i}_c \omega_X, \quad H_i^{BM}= 
 H^{-i} \omega_X.
\]
Then we have the  
following
commutative diagrams expressing the well-known
 compatibilities
of the cup and cap products: 
\be\la{tanti}
\xymatrix{
H_c^i \otimes H_c^j \ar[r] \ar[d] \ar[ddr] & H_c^{i+j} \ar[d]^= 
&&&
H_c^i \otimes H_j \ar[r] \ar[d]  \ar[ddr]& H_{j-i} \ar[d]^=
\\
H^i \otimes H_c^j \ar[r] \ar[d] \ar[dr] & H_c^{i+j}\ar[d] 
&&&
H^i \otimes H_j \ar[r] \ar[d] \ar[dr] & H_{j-i}\ar[d] 
\\
H^i \otimes H^j \ar[r]  & H^{i+j}
&&&
H^i \otimes H^{BM}_j \ar[r]  & H^{BM}_{j-i}.
}
\ee

One also has  the variants in relative cohomology
and in relative cohomology with compact supports,
the variants   with supports  on locally closed subvarieties,
as well as the variants involving a map $f: X \to Y$  (e.g. $H^*(X)$ as a $H^*(Y)$-module
etc).
The reader can sort these variants  out.

\medskip
{Author's address: Department of mathematics, Stony Brook University, 
Stony Brook, NY 11794-3651}

\end{document}